\documentclass[12pt]{article}

\newcommand{\ds}{\displaystyle}
\newtheorem{thm}{Theorem}[section]

\newcommand{\be}{\begin{equation}}
\newcommand{\ee}{\end{equation}}

\newcommand{\beq}{\begin{eqnarray}}
\newcommand{\eeq}{\end{eqnarray}}
\newcommand{\nbeq}{\begin{eqnarray*}}
\newcommand{\neeq}{\end{eqnarray*}}

\def\P{\mathbf{P}}
\def\E{\mathbf{E}}

\begin{document}

\title{Critical Branching Regenerative Processes with Migration}

\author{George P. Yanev \\
{\small Department of Mathematics} \\
{\small University of South Florida} \\
{\small St. Petersburg, FL, U.S.A.}
\and Kosto V. Mitov \\
{\small Department of Informatics and Mathematics} \\
{\small Air Force Academy 
G. Benkovski } \\
{\small Pleven, Bulgaria} \and Nick M. Yanev  \\
{\small Institute of Mathematics and Informatics} \\
{\small Bulgarian Academy of Sciences} \\
{\small Sofia, Bulgaria} }

\date { }
\maketitle 


\vspace{0.7cm} \centerline {\bf  SUMMARY}

\vspace{0.7cm}This paper demonstrates a new regeneration processes
technology making use of positive stable distributions. We study
the asymptotic behavior of branching processes with a randomly
controlled migration component. Using the new method, we confirm
some known results and establish new limit theorems that hold in a
more general setting.


\vspace{0.5cm} \noindent {\bf Key Words and Phrases}: stable laws;
alternating regenerative processes; limiting distributions;
branching processes with random migration.

\noindent {\bf AMS 2000 Subject Classification}: Primary 60K05,
60J80, Secondary 60F05.

\setcounter{equation}{0}

\setcounter{equation}{0}
\section{Introduction} \

\vspace{-0.5cm}It is well-known that stable distributions play an
important role in probability theory in general and  in the theory
of summation of random variables in particular. Zolotarev (1983)
and Uchaikin and Zolotorev (1999) give a comprehensive review of
the available results in this research area. As it is pointed out
in these references, stable laws have an extremely rich and
diverse set of applications in stochastic modelling. This holds
especially for the class of positive stable laws with parameter $0
< \alpha \le 1$ when the mean is infinite.

Mitov (1999) and Mitov and Yanev (2001) extended some of the
classical results for alternating regenerative processes that
involve positive stable distributions. Here, applying this new
regenerative technology, we obtain limit results for branching
processes with migration.

Section 2 introduces an alternating regenerative process, which
can be described as follows. The process stays at zero random
time, called ''down-period'', then it jumps up to a positive level
and re-enters the state zero after random time, called
''up-period'' or ''lifetime''. Thus, the renewal time structure of
the process consists of two components: "down" and "up", which
constitute a regeneration period. The process regenerates itself
over consecutive regeneration periods with independent and
identically distributed replicas.

The above construction applies to a variety of stochastic
processes. One example is the class of branching processes with
state dependent immigration. Mitov and Yanev (2002), using
regenerative methods, obtained limit theorems for the complex
model of Bellman-Harris branching processes with state dependent
immigration and infinite offspring variance. Another application
occurs in the area of branching diffusion processes, see Li
(2000).

In Section 3 we consider branching processes with random migration
and discuss some known results proved using traditional analytical
branching theory's methods. In Section 4  we put the migration
processes in a more general setting and extend the results from
Section 3 applying regeneration techniques. A further extension is
considered in Section 5.

Besides being of independent interest, the presented results
provide yet another example of models where the new regenerative
methods apply successfully. We believe that the scope of possible
applications is by no means limited to the listed above classes of
processes.

\setcounter{equation}{0}
\section{Alternating Regenerative Processes}

Following the description in Wolff (1989), let us have a
replacement model in which the replacement is not instantaneous.
Namely, consider a machine that breaks down and is repaired. Let
$\{T_{u,j}: \ j=1,2,\ldots \}$ be the sequence of operating times
during which the machine is ''up'' prior to a breakdown, and $\{
T_{d,j}: \ j=1,2,\ldots\}$ be the sequence of lag periods during
which the machine is ''down'' prior to completion of replacement
or repair. Assume that these sequences are independent of each
other, and that the random variables in each sequence are i.i.d.
Define $\{T_j\}$ by $T_j=T_{d,j}+T_{u,j}, \ j=1,2,\ldots$ Thus, if
the first up-period begins at $T_{d,1}$, then the first breakdown
occurs at time $T_1=T_{d,1}+T_{u,1}$. After the replacement (or
repair) the machine will start working again at $T_1+T_{d,2}$
until  the end of the second up-period $T_2=T_1+T_{d,2}+T_{u,2}$,
and so on. Call $T_j$ a repair cycle, and consider the renewal
process $N(t)$ generated by the sequence of times between
successive replacements $\{T_j\}$, i.e., \be \label{N_t}
N(t)=\max\{n\geq 0: S_n\leq t\}, \ee
 where \[S_0=0, \qquad S_n=\sum_{j=1}^nT_j,
\qquad  n=1,2,\ldots
\]
$N(t)$ is the number of repairs (replacements) completed by time
$t$.

Let us associate with  each $T_{u,j}$ a stochastic process
$\{z_j(t):0\leq t\leq T_{u,j}\} $, called cycle (or tour), $j=1,2,
\ldots$, such that
\[
z_j(0)\geq 0, \ \ z_j(t)>0 \  \mbox{for} \ 0<t<T_{u,j}, \ \
z_j(T_{u,j})=0.
\]
Assume that each $z_j(t)$ has state space $(R^+,\mathcal{B}^+)$
where $R^+=[0,\infty)$ and $\mathcal{B}^+$ is the Borel
$\sigma$field. The cycles are mutually independent and
stochastically equivalent. Also, for each $j$, \ $z_j(t)$ may
depend on $T_{u,j}$ but is independent of $\{T_{u,i}:i\neq j\}$.

Further on, consider the process $\{\sigma(t)\}$ defined by \be
\label{sigma} \sigma(t)=t-S_{N(t)}-T_{d,N(t)+1}. \ee
 Clearly,
$\sigma (t)$ can be positive or negative and
$\sigma(t)=\sigma^+(t)-\sigma^-(t)$, where
$\sigma^+(t)=\max\{\sigma(t),0\}$ is the attained duration of the
up-period in progress, whereas $\sigma^-(t)=\max\{-\sigma(t),0\}$
is the remaining time till the end of the down-period in progress.

Now, we are in a position to define an alternating regenerative
process $\{Z(t)\}$.

{\bf Definition 1}\  An alternating regenerative process {\it $\{
Z(t): t\geq 0\}$ is defined by
\[
Z(t)=
    \left\{
    \begin{array}{ll}
z_{N(t)+1}(\sigma(t)) &  \mbox{when $\sigma(t) \geq 0$}, \\
0 & \mbox {when $\sigma(t) <0$.}
    \end{array}
        \right.
\]
}

An example of alternating regenerative process is the regenerative
process with a reward structure discussed in Wolff (1989), Chapter
2. Another example is provided by the Bellman-Harris branching
processes with immigration at zero only, studied by Mitov and
Yanev (2002).

Further on, we will need three groups of assumptions.

{\bf Assumptions A} \ For the down-time component $T_{d,j}$ with
cdf $A(x)$ we assume either \[ \E T_{d,j}<\infty \]
 or
 \be \label{alpha} \E T_{d,j}=\infty \ \ \mbox{and} \ \
1-A(t)\sim t^{-\alpha}L_A(t) \ \ \mbox{as $t\to \infty$}, \ee
 where $ \alpha
\in (\frac{1}{2},1]$ and  $L_{A}(\cdot)$ is a slowly varying
function at infinity (svf), i.e., $T_{d,j}, j\geq 1$ belong to the
normal domain of attraction of a stable law with parameter
$\alpha$.

{\bf Assumptions B} \ For the up-time component $T_{u,j}$ with cdf
$F(x)$ we assume either \[ \E T_{u,j}<\infty \]
 or
 \be \label{beta}\E T_{u,j}=\infty \ \ \mbox{and} \ \
1-F(t)\sim t^{-\beta}L_F(t) \ \ \mbox{as $t\to \infty$}, \ee
 where $ \beta
\in (\frac{1}{2},1]$ and  $L_{F}(\cdot)$ is a svf, i.e., $T_{u,j},
j\geq 1$ belong to the normal domain of attraction of a stable law
with parameter $\beta$.

{\bf Assumptions C} \ For the  cycle $\{z_j(t):0\leq t\leq
T_{u,j}\} $, where $j=1,2,\ldots$ we assume \be\label{cycle}
\lim_{t \to \infty}\mathbf{P}\left\{\frac{z_j(t)}{R(t)} \le
x|T_{u,j}
>t\right\} = D(x), \ee where $R(t)=L(t)t^\gamma$, $\gamma \ge 0$ for some svf
$L(t)$ and $D(x)$ is a proper cdf on $(0,\infty)$.

\vspace{0.5cm}\noindent {\bf Basic Regeneration Theorem (BRT)}
{\it (Mitov and Yanev (2001), Mitov (1999)) Let Assumptions A -- C
hold. Set
\[
c= \lim_{t\to \infty}\frac{1-A(t)}{1-F(t)}.
\]

{\bf I.}\ Suppose that $\E T_{u,j}$ is infinite and (\ref{beta})
holds with $\frac{1}{2} < \beta < 1$.  Let $x\ge 0$.

{\bf a.}\ If $0 \le c <\infty$, then  \be\label{3.3} \lim_{t \to
\infty} \mathbf{P}\left\{\frac{Z(t)}{R(t)} \le x\right\} =
\frac{c+G(x)}{c+1}, \ee where \be\label{3.4}
G(x)=\frac{1}{B(1-\beta, \beta)} \int_0^1
D(xu^{-\gamma})u^{-\beta}(1-u)^{\beta-1}du, \ee and $B(\cdot ,
\cdot )$ stands for Beta function.

{\bf b.}\ If $c=\infty$, then  \be\label{3.5} \lim_{t \to \infty}
\mathbf{P}\left\{\frac{Z(t)}{R(t)} \le x|Z(t)>0\right\} =
\frac{1}{B(1-\beta,\alpha)} \int_0^1
D(xu^{-\gamma})u^{-\beta}(1-u)^{\alpha-1}du. \ee

{\bf II.}\ Suppose that $\E T_{u,j}$ is infinite and (\ref{beta})
holds with $\beta=1$. Assume (\ref{cycle}) with $D(0)=0$ and let
$0<x<1$.

{\bf a.}\ If $0 \le c < \infty$, then  \be\label{3.7} \lim_{t \to
\infty}\mathbf{P}\left\{\frac{m_F(R^{-1}(Z(t)))}{m_F(t)} \le
x\right\} = \frac{c+x}{c+1}, \ee where $R^{-1}(\cdot)$ is the
inverse function of $R(\cdot)$ and $m_F(t)=\int_0^t1-F(x)dx$.

{\bf b.}\ If $c = \infty$, then \be\label{3.8} \lim_{t \to
\infty}\mathbf{P}\left\{\frac{m_F(R^{-1}(Z(t)))}{m_F(t)} \le
x|Z(t)>0\right\} = x, \ee where $R^{-1}(\cdot)$ and $m_F(t)$ are
defined above in part {\bf a}.

 {\bf III.}\ Suppose that $\E T_{d,j}$ is infinite and
(\ref{alpha}) holds. If $\E T_{u,j}<\infty$, then \be\label{3.9}
\lim_{t \to \infty} \mathbf{P}\left\{Z(t) \le x|Z(t)
>0\right\} = \frac{1}{\E T_{u,1}}\int_0^\infty \mathbf{P}\left\{z_1(y) \le x,
T_{u,1}>y\right\} dy. \ee }

{\bf Remark}\  Notice that, if both $\E T_{d,j}$ and $\E T_{u,j}$
are finite, then by the classical regeneration theorem (see e.g.
Sigman and Wolf (1993), Theorem 2.1) \be\label{3.10} \lim_{t \to
\infty} \mathbf{P}\left\{Z(t) \le x\right\} = \frac{1}{\E
T_{d,1}+\E T_{u,1}}\int_0^\infty
\mathbf{P}\left\{z_1(y) \le x,T_{d,1}+T_{u,1}>y\right\} dy. \ee  

The BRT in the non-lattice case was proved by Mitov and Yanev
(2001) and in the lattice case by Mitov (1999). In the next
section, applying the BRT, we obtain limit theorems for a class of
branching processes with random migration.

\setcounter{equation}{0}
\section{Branching Processes with Migration}

Let us have on a probability space $( \Omega ,{\cal F} , P)$ three
independent sets of non-negative integer-valued i.i.d. random
variables as follows.

{\bf i.}\  offspring variables: $\{ X_{it}: \ i=1,2,\ldots;
t=0,1,\ldots\}$;

 {\bf ii.}\ immigration variables: $\{ (I_t^+, I_t^o): \
t=0,1,\ldots\}$;

 {\bf iii.}\ emigration variables: $\{ (_{fam}\mathcal{E}_t, \  _{ind}\mathcal{E}_t): \
t=0,1,\ldots \}$.

Let us construct a sequence $\{ (M_t^+,M_t^o): \ t=0,1,\ldots \}$,
which will play the role of a migration component of the process.
Set
\[ M_t^+ = \left \{
\begin{array}{ll}
       {\ds -\sum_{i=1}^{\ds _{fam}\mathcal{E}_t}X_{it}- \ _{ind}\mathcal{E}_t } & \mbox{with probability
\ {\it p},} \\
           0                    & \mbox{with probability \ {\it q},} \\
I_t^+       & \mbox{with probability \  {\it r}, \qquad \qquad {\it p}+{\it q}+{\it
r}=1}
            \end{array}
     \right. \]
and
\[ M_t^o = \left \{
\begin{array}{ll}

           0                    & \mbox{with probability \
1-{\it r},} \\
I_t^o       & \mbox{with probability \  {\it r}.}
            \end{array}
     \right. \]

Let us define a branching process with migration.

{\bf Definition 2}\  {\it A branching process with randomly
controlled migration $\{ Y_t: \ t=0,1,\ldots \}$ is defined by the
recurrence \be \label{migration} Y_{t+1}=\max\left\{
\sum_{i=1}^{Y_t}X_{it}+M_t, 0\right\},\qquad t=0,1,2,..., \ee
 where
\[
 M_t = M_t^+1_{\ds \{Y_t>0\}}+ M_t^o1_{\ds \{Y_t=0\}}
\]
and \ $Y_0\geq 0$ \ is independent of \ $X_{\cdot t}$ and $M_t$
for $t>0$. The equalities above hold in distribution and $1_A$
stands for the indicator of an event $A$. }

The process $\{ Y_t\}$ is a homogeneous Markov chain which admits
the following interpretation.  Upon the reproduction in the $t$th
generation three situations are possible:  {\bf (i)} with
probability $q$ the process develops like a
Bienaym\'{e}-Galton-Watson process, i.e., without any migration;
{\bf (ii)} with probability $p$ there is an emigration of ${\ds
_{fam}\mathcal{E}_t}$ families and  $ _{ind}\mathcal{E}_t$
individuals; {\bf (iii)} with probability $r$ there is an
immigration of $I_t^+$ or $I_t^o$ new particles depending on the
state of the process.


Discrete time processes with different migration components were
introduced by Yanev and Mitov (1980, 1985) and Nagaev and Han
(1980), see Rahimov (1995) for throughout discussions and
additional references. Continuous time branching processes
regulated by different schemes of emigration ("catastrophes",
"disasters") and immigration have also been subject of
considerable interest. Let us point out here Pakes (1986), Chen
and Rensaw (1995), and Rahimov and Al-Sabah (2000) papers among
others.

Note that, definition (\ref{migration}) includes as its particular
cases some well-known models: processes with immigration (when
$r=1$ and $I_t^+\equiv I_t^o$ ), processes with state-dependent
immigration (when $r=1$ and $I_t^+\equiv 0$), and processes with
pure emigration (when $p=1$).

We shall study branching processes with migration assuming
offspring mean one (critical case) and finite variance, i.e.,
\be\label{critical} \E X_{it}=1 \quad \mbox{and} \quad
0<VarX_{it}=2b <\infty , \quad \mbox{say}. \ee
 Previous studies
(see Yanev and Yanev (1995 - 1997)) revealed the importance of a
parameter that relates both reproduction and migration components,
given by
\[
\theta= \frac{\E M_t^+}{b}.
\]
Note that, a similar parameter that measures the relative sizes of
immigration and branching in case of processes with immigration,
first appeared in Zubkov~(1972). $\theta$ is the recurrence
parameter of the Markov chain: $\{ Y_t\}$ is non-recurrent for
$\theta>1$ (when immigration strongly dominates emigration); it is
null-recurrent for $0\leq \theta < 1$ (when immigration mildly
dominates emigration), and positive recurrent for $\theta <0$
(when emigration dominates immigration). In the border case
$\theta =1$ the chain is either non-recurrent or null-recurrent
depending on some extra moment assumptions.

Further on, we need additional assumptions for the migration
component as follows \be\label{ass_migration}
\begin{array}{l}
0 < \mathbf{E}I_t^+ <\infty, \ \ 0<\mathbf{E}I_t^o <\infty,
\\
0 \le \ _{fam}\mathcal{E}_t \le C_1 < \infty, \ \ \ 0 \le \
_{ind}\mathcal{E}_t \le C_2 <\infty, \ \mbox{ a.s.}.
\end{array}
 \ee

The following theorem gives the limiting behavior of $\{ Y_t\}$.

\begin{thm}\label{old} (Yanev and Yanev (1996)) Suppose that $\{ Y_t\}$ is critical with finite offspring
variance, i.e., (\ref{critical}). Also, assume that the migration
satisfies (\ref{ass_migration}).

{\bf I.} If $\theta>0$, then
$$ \frac{Y_t}{bt} \stackrel{d}{\to} \Gamma(\theta,1),$$
where the limit is Gamma distributed with parameters $\theta$ and
$1$.

{\bf II.} If $\theta=0$ and $\E I_t^{+2}<\infty$, then
$$\frac{\log Y_t}{\log t} \stackrel{d}{\to } U(0,1), $$ where the
limit is uniformly distributed on the unit interval.

{\bf III.} If $\theta <0$, then $Y_t$ possesses a limiting
stationary distribution, i.e.,
$$Y_t \stackrel{d}{\to} Y_\infty .$$
\end{thm}

The above theorem was proved using some traditional branching
process theory methods including functional equations for pgf's,
Laplace transforms, Tauberian theorems etc. In the next section we
will extend the above results making use of some probabilistic
arguments and the BRT from Section 2.

\setcounter{equation}{0}
\section{Regeneration and Migration}

Consider the branching migration process  (\ref{migration}) with
$M_t^o\equiv 0$, i.e., migration is not permitted when the process
is in state zero. That is, \be \label{stopped}
Y^o_{t+1}=\max\left\{
\sum_{i=1}^{Y^o_t}X_{it}+M^+_t1_{\{Y_t^0>0\}},0\right\},\qquad
t=0,1,2,..., \ee where \ $Y^o_0\geq 0$ \ is independent of \
$X_{\cdot t}$ and $M^+_t$ for $t\geq 1$. Call this a process with
migration stopped at zero.

Since the state zero is a reflective barrier, the Markov chain
$\{Y_t\}$ is an regenerative process. Indeed, using the notations
from Section 2, $\{ Y_t\}$ stays at zero random time $T_{d,1}$,
which has the geometric distribution \be \label{geo}
\mathbf{P}\{T_{d,1}=k\}=\P^{k-1}\{M_t^o=0\}(1-\P\{M_t^o=0\}), \ \
k=1,2,\ldots \ee
 In the end of this down-period the process jumps
up to a random level $I_{T_{d,1}}^0$ and evolves according to the
rules in the model (\ref{migration}) until it hits zero again in
the end of its lifetime $T_{u,1}$. Thus, $T_1=T_{d,1}+T_{u,1}$
forms the first period of regeneration and the evolution of the
process repeats in the next such periods, i.e., $\{ Y_{t+T_1}: \
t\geq 0\}$ is stochastically equivalent to $\{ Y_{t}: \ t\geq
0\}$. Let $T_{d,j}$, $j=1,2,\ldots$ be i.i.d. copies of $T_{d,1}$
given by (\ref{geo}). Also, let $\{ Y_{j,t}^o: j=1,2,\ldots \}$ be
a sequence of branching processes with migration stopped at zero
defined by (\ref{stopped}), having lifetimes $T_{u,j}$, i.e., for
$j=1,2,\ldots$
\[
Y_{j,0}^o\geq 0,  \ \ Y_{j,t}^o>0 \  \mbox{for} \ 0<t<T_{u,j}, \ \
Y_{j,T_{u,j}}^o=0.
\]
Now, it is not difficult to see, that (\ref{migration}) is a
particular case of Definition 1 with cycle process $\{ Y_{j,t}^o:
t=0,1,\ldots T_{u,j} \}$. Indeed, $\{Y_t\}$ is an alternative
regenerative process with cycle process $z_j(t)\equiv Y_{j,t}^o$
and $T_{d,j}$ given by (\ref{geo}).

Let us generalize the migration process (\ref{migration}) as
follows:

{\bf i.}\ First, assume that the down-periods $\{ T_{d,j}: \
j=1,2,\ldots\}$ have a cdf $A(x)$, which is not necessarily the
geometric one from (\ref{geo}).

{\bf ii.}\ Secondly, assume that the mean of $\{ T_{d,j}: \
j=1,2,\ldots\}$ is not constrained to be finite. If $\E
T_{d,j}=\infty$, assume that $T_{d,j}$'s belong to the domain of
attraction of a stable law and their cdf $F(t)$ satisfies
(\ref{beta}).

Appealing to Definition 1, we construct a generalized version of
$\{ Y_t\}$ with {\bf i.} and {\bf ii.} above as follows

{\bf Definition 3}\ {\it A  branching regenerative process with
migration denoted by \newline $\{ Z_t: t=0,1,\ldots \}$ is defined
by
\[
Z_t=
    \left\{
    \begin{array}{ll}
Y^o_{N(t)+1,\ \sigma(t)} &  \mbox{when $\sigma(t) \geq 0$}, \\
0 & \mbox {when $\sigma(t) <0$,}
    \end{array}
        \right.
\]
where the cycles $\{ Y^o_{j,t}\}$ are processes with migration
stopped at zero;
 $N(t)$ and $\sigma(t)$ are defined by
(\ref{N_t}) and (\ref{sigma}), respectively. }

{\bf Example}\ Using the well-known duality between a branching
process and a M/M/1 queue, let us describe a situation where the
above construction applies. Consider a queueing model in which
customers arrive following a Poisson process. Then the successive
times $T_j$ from the commencement of the $j$th busy period to the
start of the next busy period form a renewal process. Each $T_j$
is composed of a busy portion $T_{u,j}$ and an idle portion
$T_{d,j}$, when the queue is not empty or empty, respectively.
Assuming that the customers, arriving during the service time of a
customer, are his/her "offspring", we obtain a branching
regenerative process. The immigration component accounts for a
policy when certain customers (probably coming from a second
source) accumulate and will be served after completing the service
time of a "generation". Alternatively, some customers (called
"emigrants") may leave the system prior to their service.

Further on, we will assume that some reproduction and immigration
moments are finite as follows

\be\label{moments}
\begin{array}{lll}
  \mathbf{E}I_t^{+2} < \infty , & \ & \mbox{when}\
\theta  =0 \\
\E I_t^{+(1-\theta)}<\infty  , & \E X_{1t}^2\log (1+X_{1t})<\infty
& \mbox {when} \  -1 < \theta
< 0  \\
\E I_t^{+2}\log^2(1+I_t^{+})<\infty , & \E
X_{1t}^2\log^2(1+X_{1t})<\infty
& \mbox {when} \   \theta  = -1    \\
 \E I_t^{+(1-\theta)}<\infty , & \E X_{1t}^{1-\theta}<\infty    &
\mbox {when} \   \theta < -1 .
   \end{array}
\ee

Let us summarize for more convenient references some results for
processes with migration stopped at zero 
that are
proved in Yanev and Yanev~(1995-2002).

\begin{thm}\label{ThmCycle} (Yanev and Yanev (1995-2002))
Suppose that $\{ Y_t^o\}$ is critical with finite offspring
variance, i.e., (\ref{critical}). Also, assume that the
reproduction and migration satisfy (\ref{ass_migration}) and
(\ref{moments}).

If $\theta\geq 0$, then \be\label{T_u}
 \P\{T_{u,1}>0\}\sim L(t)t^{-(1-\theta)\vee 0},
\ee where $L(t)$ is a svf and

\be\label{Y^o}
 \frac{Y^o_t}{bt}\ |\ Y^o_t>0
\stackrel{d}{\to} Y^o_\infty, \ee
 where the limiting random
variable $Y^o_\infty$ has Exp(1) distribution when $-\infty <
\theta \leq 1$ and Gamma $(\theta,1)$ when $\theta>1$.
\end{thm}

Therefore, in effect, if the immigration's domination is
insufficient ($\theta \leq 1$) to prevent certain extinction, the
conditioned Kolmogorov-Yaglom's exponential limit law for
processes without any migration holds; whereas if sufficient, the
limit law coincides with that in Theorem~\ref{old}i, where
immigration is also permitted in zero, on the set of
non-extinction.

Recall, from the BRT in Section 2, that
\[
c=\lim_{t\to \infty}\frac{\P\{T_{d,j}>0\}}{\P\{ T_{u,j}>0\}}.
\]
Assuming that at least one of $\E T_{d,j}$ and $\E T_{u,j}$ is not
finite, one can interpret the parameter $c$ as follows. If $0\leq
c<\infty$, then the up-period (process' lifetime) has length
asymptotically bigger than the down-period (stay at zero); whereas
if $c=\infty$, then the process spends more time at zero then up.

Now, we are in a position to prove the main limit theorem for
branching regenerative processes with migration.

\begin{thm}\label{main}
 Suppose that $\{
Z_t\}$ is critical with finite offspring variance, i.e.,
(\ref{critical}). Also, assume that the reproduction and migration
satisfy (\ref{ass_migration}) and (\ref{moments}).

{\bf I.} Let $0< \theta < 1/2$ and suppose that either $\E
T_{d,j}<\infty$ or $\E T_{d,j}=\infty$ and (\ref{alpha}) holds.

{\bf a.}\  If $0\leq c <\infty$, then

\be\label{0theta1}
 \frac{Z_t}{bt} \stackrel{d}{\to} Z_\infty,
\ee
 where $\mathbf{E}Z_\infty=\theta/(c+1)$ and for $x \geq 0$
\be\label{0theta1limit}
 \P \{ Z_\infty \leq x\} =
1-\frac{1}{(c+1)B(\theta,1-\theta)}\int_0^1e^{-x/y}y^{\theta-1}(1-y)^{-\theta}
dy. \ee

{\bf b.}\  If $c =\infty$, then

\be\label{0theta12}
 \frac{Z_t}{bt} \ | \ Z_t>0 \stackrel{d}{\to} Z_\infty,
\ee
 where $\mathbf{E}Z_\infty=\theta/(\theta + \alpha)$ and for $x \geq 0$
\be\label{0theta1limit2}
 \P \{ Z_\infty \leq x\} =
1-\frac{1}{B(\theta,\alpha)}\int_0^1e^{-x/y}y^{\theta-1}(1-y)^{\alpha-1}
dy .\ee

{\bf II.} Let $\theta=0$.

{\bf a.}\  If \ $\E T_{d,j}<\infty$, then
$$\frac{\log
Z_t}{\log t} \stackrel{d}{\to } U(0,1), $$ where the limit is
uniformly distributed on the unit interval.

{\bf b.}\  Assume $\E T_{d,j}=\infty$ and (\ref{alpha}). {\bf i.}
If $ 0\leq c<\infty$, then
$$\frac{\log
Z_t}{\log t} \stackrel{d}{\to } Z_\infty, $$ where $\P \{ Z_\infty
\leq x\}= (c+x)/(c+1)$ for $0\leq x\leq 1$.

\noindent {\bf ii.} If $ c=\infty$, then
$$\frac{\log
Z_t}{\log t} \ | \ Z_t>0 \ \stackrel{d}{\to } U(0,1), $$ where the
limit is uniformly distributed on the unit interval.

 {\bf III.} Let $\theta < 0$.

 {\bf a.}\  If \ $\E T_{d,j}<\infty$, then $Z_t$ possesses a limiting
stationary distribution, i.e.,
$$Z_t \stackrel{d}{\to} Z_\infty $$
and for $x\geq 0$
\[
\mathbf{P}\{Z_\infty \leq x\}= \frac{1}{\E T_{1}}
\sum_{k=0}^\infty \mathbf{P}\{Y^o_{1,k}\leq x , \  T_{1}> k\} .
\]

{\bf b.}\  If \ $\E T_{d,j}=\infty$ and (\ref{alpha}) holds, then

$$Z_t \ | \ Z_t>0 \stackrel{d}{\to} Z_\infty $$
and for $x\geq 0$
\[
\mathbf{P}\{Z_\infty \leq x\}= \frac{1}{\E T_{u,1}}
\sum_{k=0}^\infty \mathbf{P}\{Y^o_{1,k}\leq x , \  T_{u,1}> k\} .
\]
\end{thm}

{\bf Proof}\ We shall apply the BRT to $\{ Z_t\}$. {\bf I.} First
note that, (\ref{T_u}) implies \[ \mathbf{P}\{T_{u,1}
> t \} \sim L(t)t^{-(1-\theta)}
\] and hence (\ref{beta}) holds with $\beta=1-\theta$. Furthermore,
by (\ref{Y^o}), \be\label{4.22} \lim_{t \to
\infty}\mathbf{P}\left\{\frac{Y^o_t}{bt} \le x|Y^o_t >0\right\}
=1-e^{-x}. \ee Thus, (\ref{cycle}) holds with $R(t)=bt$ and
$D(x)=1-e^{-x}$. Now, (\ref{3.3}) and (\ref{3.4}) lead to
(\ref{0theta1}) with \nbeq \P\{Z_\infty\leq x\} & = &
\frac{c}{c+1}+ \frac{1}{(c+1)B(\theta,
1-\theta)}\int_0^1(1-e^{x/u})u^{-(1-\theta)}(1-u)^{-\theta}du \\
    & = &
 1-   \frac{1}{(c+1)B(\theta,1-\theta)}\int_0^1e^{x/u}u^{\theta-1}(1-u)^{-\theta}du,
\neeq which proves (\ref{0theta1limit}). Integrating
$1-\P\{Z_\infty\leq x\}$, it is not difficult to obtain
$EZ_\infty$. Similarly, taking into account (\ref{4.22}) and
(\ref{3.5}), one can obtain (\ref{0theta12}) and
(\ref{0theta1limit2}).

{\bf II.}\ In this case we have from (\ref{T_u}) that
\[
\mathbf{P}\{T_{u,1}
> t \} \sim C_0 t^{-1},
\]
 for some positive constant $C_0$. Thus, (\ref{beta}) holds with
$\beta=1$ and $m_f(t)\sim C_0\log t$. On the other hand,
(\ref{4.22}) still applies and hence $R(t)=bt$ and $D(x)=1-e^{-x}$
with $D(0)=0$ . Finally,
\[
\frac{m_F(M^{-1}(Z_t))}{m_F(t)} \sim \frac{\log Z_t}{\log t}.
\]
Now, part {\bf IIa} follows from (\ref{3.7}), taking into account
that $\E T_{d,j}<\infty$, which results in $c=0$. Similarly,
(\ref{3.7}) and (\ref{3.8}) imply {\bf IIbi} and {\bf IIbii},
respectively.

{\bf III.}\ According to Theorem~\ref{ThmCycle}, we have $\E
T_{u,j}<\infty$. In case of $\E T_{d,j}<\infty$, the classical
regeneration theorem applies and (\ref{3.10}) leads to {\bf IIIa}.
If $\E T_{d,j}=\infty$, then (\ref{3.9}) holds and hence {\bf
IIIb}. \hfill{\rule{1.6ex}{1.6ex}}

It is interesting to compare  Theorem~\ref{old} for $\{ Y_t\}$
with Theorem~\ref{main} for $\{ Z_t\}$. Recall that the former
assumes a geometric distributed $T_{d,j}$ with $\E
T_{d,j}<\infty$, whereas the later holds for $T_{d,j}$, not
constrained to one specific distribution and that might have a
finite or infinite mean. For all values of $\theta$, if $\E
T_{d,j}=\infty$ and $c=\infty$, i.e., asymptotically the
down-period dominates the up-period, Theorem~\ref{main} represents
new conditional limit results, on the set of non-extinction. The
case of finite $\E T_{d,j}$ results in unconditional limit results
and if $\theta \leq 0$ then  the results in Theorem~\ref{old}II
extend to the more general process $\{ Z_t\}$. In the intermediate
case when $\E T_{d,j}=\infty$ but $0\leq c<\infty$, i.e., the
up-period dominates the down one, we obtain unconditional limit
results when $\theta \geq 0$ and a conditional one when
$\theta<0$.

It is worth mentioning the generality we gain in
Theorem~\ref{main} due to the new regeneration methods of proof
versus the traditional probability generation functions based
techniques used to obtain Theorem~\ref{old}. One limitation of the
new approach is that it does not apply for $\theta \geq 1/2$.

\setcounter{equation}{0}
\section{One Extension}

The results presented in Section 4 can be extended by relaxing the
condition that $\E Y_{j,t}^o<\infty$. Instead, let us assume that
the immigration at zero belongs to the domain of attraction of a
stable law with parameter $0<\rho\leq 1$. For a similar extension
in case of processes with immigration at zero only, see Ivanoff
and Seneta~(1985). Further on, we will need the following result.

\begin{thm}\label{ThmCycleRho} (Yanev and Yanev (1997))
Suppose that $\{ Y_t^o\}$ is critical with finite offspring
variance, i.e., (\ref{critical}). Also, assume that the
reproduction and migration satisfy (\ref{ass_migration}) and $\E
I_t^{+2}<\infty$ when $\theta=0$.

{\bf I.}\ If $\theta + \rho \geq 1$, then
\be\label{T_uRho1}
 \P\{T_{u,1}>0\}\sim K(t)t^{-(1-\theta)\vee 0},
\ee where $K(t)$ is a svf and \[
 \frac{Y^o_t}{bt}\ |\ Y^o_t>0
\stackrel{d}{\to} Y^o_\infty, \]
 where the limit $Y^o_\infty$ has Exp(1) distribution when $-\infty < \theta \leq
1$ and Gamma $(\theta,1)$ when $\theta>1$.

{\bf II.}\ If $\theta + \rho < 1$, then \[
 \P\{T_{u,1}>0\}\sim K(t)t^{-\rho},
\] where $K(t)$ is a svf and \be\label{Y^oRho2}
 \frac{Y^o_t}{bt}\ |\ Y^o_t>0
\stackrel{d}{\to} Y^o_\infty, \ee
 where the limit
$Y^o_\infty$ has a proper cdf with Laplace transform given by
\be\label{varphi} \varphi(\lambda)=1-\frac{ \lambda^\rho (1-\theta
-\rho)}{(1+\lambda)^{\theta+ \rho}}B(1-\rho, 1-\theta) -\lambda
\theta \int_0^1\frac{1}{(1-y)^{\rho}(1+\lambda y)^{\theta +1}}dy.
\ee

\end{thm}

Let us point out that, if $\theta+\rho<1$ then the rate of
convergence in (\ref{T_uRho1}) depends on $\rho$ only. One can say
that $I_t^o$, the ancestors' distribution for $Y_t^o$, dominates
the migration component. Indeed, in this case the limit in
(\ref{Y^oRho2}) and $I_t^o$ share the domain of attraction of the
same stable law with parameter $\rho$. If $\theta + \rho \geq 1$,
then the form of the limiting distribution depends essentially on
$\theta$ and one might say that in this case the migration
component is dominating.

Applying the BRT and Theorem~\ref{ThmCycleRho}, we obtain the
following limit theorem.

\begin{thm}\label{mainRho} Suppose that $\{
Z_t\}$ is critical with finite offspring variance, i.e.,
(\ref{critical}). Also, assume that the reproduction and migration
satisfy (\ref{ass_migration}) and $\E I_t^{+2}<~\infty$ when
$\theta=0$.

{\bf I.}\  If  $0 < \theta < 1/2$ and $1/2<\rho<1$, such that
$\rho+\theta \ge 1$, then the limit results in Theorem~\ref{main}I
extend to $\{ Z_t\}$, under the same assumptions on $c$.

{\bf II.}\ Assume $ \theta < 1/2$ and $1/2<\rho<1$, such that
$\rho+\theta < 1$.

{\bf a.}\ If $0\leq c<\infty$, then
\[
\frac{Z_t}{bt} \stackrel{d}{\to} Z_\infty,
\]
 where the limit $Z_\infty$ has Laplace transform
\beq \label{LT1}
 \lefteqn{\E e^{-\lambda Z_\infty}}
\\
\! \! & \! \! = \! \! & \frac{1}{c+1} \left( 1 - \frac{I_{\lambda
/\lambda+1}(\rho, 1-\theta-\rho) }{(\lambda+1)^\theta } - \frac{
\lambda \theta} {B(1-\rho,\rho)}\int_0^1 \int_0^1
\frac{u^{1-\rho}(1+\lambda y u)^{-\theta-1 }}{(1-u)^{1-\rho}
(1-y)^\rho}dy du \right) \nonumber
\eeq where $I_x(a,b)=B_x(a,b)/B(a,b)$.

{\bf b.} \ If $c=\infty$, then
\[
 \frac{Z_t}{bt} \ | \ Z_t>0 \stackrel{d}{\to} Z_\infty,
 \]
 where the limit $Z_\infty$ has Laplace transform
\beq \label{LT2}
 \lefteqn{\E e^{-\lambda Z_\infty}}
\\
\! \! & \! \! = \! \! &  1 - \frac{I_{\lambda /\lambda+1}(\alpha,
1-\theta-\rho)}{C_{\lambda}(\alpha, \theta, \rho)}  - \frac{
\lambda \theta }{B(1-\rho,\alpha)}\int_0^1 \int_0^1
\frac{u^{1-\rho}(1+\lambda y u)^{-\theta-1
}}{(1-u)^{1-\alpha}(1-y)^\rho}dy du, \nonumber \eeq where
$C_{\lambda}( \alpha,
\theta,\rho)=\lambda^{\alpha-\rho}(\alpha+1-\theta-\rho)/\left((\lambda+1)^{\alpha
- \theta -\rho} B(1-\theta, \alpha + 1 - \rho)\right) $.

If either $c=0$ or $c=\infty$, then the limiting distributions in
(\ref{LT1}) and (\ref{LT2}) belong to the normal domain of
attraction of a stable law with parameter $\rho \in
(\frac{1}{2},1).$

{\bf III.}\ If $\theta \le 0$ and $\rho=1$, then the limit results
in Theorem~\ref{main}II extend to $\{ Z_t\}$, under the same
assumptions on $T_{d,j}$ and $c$ .

\end{thm}

{\bf Proof}\ {\bf I.}\ Under the hypotheses,
Theorem~\ref{ThmCycleRho} implies
$
 \mathbf{P}\{T_{u,1}
> t \} \sim K(t)t^{-(1-\theta)}
$ and $
 \lim_{t \to \infty}\mathbf{P}\left\{Y^o_t/(bt) \le x|Y^o_t
>0\right\} =1-e^{-x}.
$
The rest of the proof repeats the arguments in that of
Theorem~\ref{main}I.

{\bf II.}\ Let us apply the BRT again. According to
Theorem~\ref{ThmCycleRho}
\[
 \mathbf{P}\{T_{u,1}
> t \} \sim K(t)t^{-\rho}.
\]
Thus, (\ref{beta}) holds with $\beta=\rho$. On the other hand, by
the same theorem,
\[
 \lim_{t \to \infty}\mathbf{P}\left\{\frac{Y^o_t}{bt} \le x|Y^o_t
>0\right\} =Y^o_\infty
\]
and the limiting random variable has a proper cdf $H(x)$ with
Laplace transform given by (\ref{varphi}). Let $0\leq c<\infty$.
Then (\ref{3.3}) with (\ref{3.4}) (note that $\gamma=1$ and $\beta
=\rho)$ implies the limiting result {\bf IIa} and \nbeq \E
e^{-\lambda Z_\infty}
 & = & \int_0^\infty e^{-\lambda x}d\frac{c+G(x)}{c+1} \\
 & = & \frac{\lambda}{c+1}\int_0^\infty e^{-\lambda x}G(x)dx \\
 & = & \frac{\lambda}{(c+1)B(1-\rho,
 \rho)}\int_0^1u^{-\rho}(1-u)^{\rho-1}
 \int_0^\infty e^{-\lambda x}D(x/u)dx du \\
 & = &
\frac{\lambda}{(c+1)B(1-\rho,
 \rho)}\int_0^1u^{-\rho}(1-u)^{\rho-1}
 \int_0^\infty e^{-\lambda uy}D(y)dy du \\
& = & \frac{1}{(c+1)B(1-\rho,
 \rho)}\int_0^1u^{-\rho}(1-u)^{\rho-1}
  \varphi (\lambda u) du ,
  \neeq
  where $\varphi (\cdot)$ is given by (\ref{varphi}).
  After replacing $\varphi(\lambda u)$ with (\ref{varphi}) and using well-known
  properties of the incomplete Beta function,
  we obtain (\ref{LT1}). The case $c=\infty$ follows
similarly by (\ref{3.5}).

Let us prove that if $c=0$ then the limiting distribution with
(\ref{LT1}) belongs to the normal domain of attraction of a stable
law with parameter $\rho$. Indeed, setting $H(x)=\P \{Z_\infty
\leq x\}$, from (\ref{LT1}) since $I_x(a,b)\sim const. \ x^a$, we
obtain as $\lambda \to 0+$ \nbeq \int_0^\infty e^{-\lambda
x}(1-H(x))du & = & \frac{1-\E e^{-\lambda Z_\infty}}{\lambda} \sim
const. \lambda^{\rho-1}. \neeq Now, by Karamata's Tauberian
theorem $ \int_0^t1-H(x)dx \sim 
const. t^{1-\rho} $ as $t \to \infty$. Thus,  $ 1-H(x) \sim
const. x^{\rho}, $ which was to be proved. The statement for
$c=\infty$ follows similarly from (\ref{LT2}).

{\bf III.}\ Since $\rho=1$, Theorem~\ref{ThmCycleRho} yields
$
\mathbf{P}\{T_{u,1}
> t \} \sim K(t) t^{-1}
$ and also $ \lim_{t \to \infty}\mathbf{P}\left\{Y^o_t/(bt) \le
x|Y^o_t>0\right\} =1-e^{-x}. $ We complete the proof of {\bf III}
by repeating the arguments in the proof of Theorem~\ref{main}II.
\hfill{\rule{1.6ex}{1.6ex}}

{\bf Remarks} \ Notice the new limiting distributions that appear
in Theorem~\ref{mainRho}II, when the immigration at zero $I_t^o$
dominates the migration component
. If either $c=0$ or $c=\infty$, then the limiting distributions
belong to the same stable law's domain as the immigration at zero
$I_t^o$. We can also deduce that they
have no mass at zero. Indeed if $c=0$, then one can show (see
Yanev and Yanev (1997)) that $\lim_{\lambda \to \infty}(c+1)(1-\E
e^{-\lambda Z_\infty})=1$ and hence $\lim_{\lambda \to \infty}\E
e^{-\lambda Z_\infty}=0$. The case $c=\infty$ is similar. Finally,
it is interesting to see that if $\theta=0$, then the Laplace
transform in {\bf IIa.} simplifies to
$(1-I_{\lambda/(\lambda+1)}(\rho, 1-\rho))/(c+1).$ This extends
the result from (\ref{varphi}) when
$\varphi(\lambda)=1-(\lambda/\lambda+1)^\rho$. If $\theta=0$ and
$c=\infty$ we have $\E e^{-\lambda
Z_\infty}=1-I_{\lambda/(\lambda+1)}(\alpha,
1-\rho))/C_\lambda(\alpha, 0, \rho) $.

\vspace{0.3cm}{\bf Acknowledgments}\  This research is supported
in part by NFSI of Bulgaria, grant No. MM 1101/2001.



\end{document}